\documentclass[final]{siamltex}
\usepackage{amsmath}
\usepackage{amssymb}
\usepackage{array}
\usepackage{blkarray}
\usepackage{capt-of}
\usepackage{color}
\usepackage{colortbl}
\usepackage{fancyhdr}
\usepackage{fancyvrb}
\usepackage{float}
\usepackage[T1]{fontenc}
\usepackage{framed}
\usepackage{geometry}
\usepackage{graphicx}
\usepackage{hyperref}
\usepackage[latin1]{inputenc}
\usepackage{stmaryrd}
\usepackage{subfigure}
\usepackage{upquote}
\usepackage{url}
\usepackage[usenames,dvipsnames,svgnames,table]{xcolor}
\usepackage{xspace}
\newcommand{\nc}{\newcommand}
\nc{\dsp}{\displaystyle}
\nc{\txt}{\textstyle}
\nc{\reff}[1]{(\ref{#1})}
\nc{\mrm}[1]{\mathrm{#1}}
\nc{\udl}[1]{\underline{#1}}
\nc{\ovl}[1]{\overline{#1}}
\nc{\al}{\underline{\boldsymbol{\alpha}}}
\nc{\la}{\underline{\boldsymbol{\lambda}}}
\nc{\llbr}{\llbracket}
\nc{\rrbr}{\rrbracket}
\nc{\lbr}{\lbrack}
\nc{\rbr}{\rbrack}
\nc{\N}{\mathbb{N}}
\nc{\Z}{\mathbb{Z}}
\nc{\D}{\mathbb{D}}
\nc{\Q}{\mathbb{Q}}
\nc{\R}{\mathbb{R}}
\nc{\C}{\mathbb{C}}
\nc{\T}{\mathbb{T}}
\nc{\Stwo}{\mathbb{S}^2}
\nc{\tld}[1]{\tilde{#1}}
\nc{\wtld}[1]{\widetilde{#1}}
\nc{\hu}{\hat{u}}
\nc{\wh}[1]{\widehat{#1}}
\nc{\sumeven}{\sum_{k=-N/2}^{N/2}{\hspace{-0.3cm}}'{\;\,}}
\nc{\sumodd}{\sum_{k=-\frac{N-1}{2}}^{\frac{N-1}{2}}}
\nc{\sumoddl}{\sum_{l=-\frac{N-1}{2}}^{\frac{N-1}{2}}}
\nc{\cqfd}{~\hbox{\vrule width 2.5pt depth 2.5 pt height 3.5 pt}}
\title{Computing planar and spherical choreographies}
\author{Hadrien Montanelli\thanks{Oxford University Mathematical Institute, Oxford OX2 6GG, UK.\ \ Supported by 
the European Research Council under the European Union's Seventh Framework Programme (FP7/2007--2013)/ERC grant agreement 
no.\ 291068.\ \ The views expressed in this article are not those of the ERC or the European Commission, and the European Union is not 
liable for any use that may be made of the information contained here.} \and Nikola I. Gushterov\thanks{Oxford University Centre for Theoretical Physics, Oxford OX1 3NP, UK.}}

\begin{document}

\maketitle

\begin{abstract}
An algorithm is presented for numerical computation of choreographies in the plane in a Newtonian potential and on the sphere in a
cotangent potential.
It is based on stereographic projection, approximation by trigonometric polynomials, and quasi-Newton and Newton optimization methods
with exact gradient and exact Hessian matrix. New choreographies on the sphere are presented. 
\end{abstract}


\pagestyle{myheadings}
\thispagestyle{plain}

\markboth{COMPUTING PLANAR AND SPHERICAL CHOREOGRAPHIES}{MONTANELLI AND GUSHTEROV}

\section{Introduction} 

Choreographies are periodic solutions of the $n$-body problem, $n\geq2$, in which the bodies share a common orbit and are uniformly spread along it. 
The study of choreographies with unit masses in the plane in a Newtonian potential is an old one. 
For the two-body problem, the only choreography is a circle, and for the three-body problem, the first one was found by Lagrange in 1772~\cite{lagrange1772}, and is also a circle.
The second choreography of the three-body problem with unit masses, a figure-eight, was discovered numerically more than two centuries later by Moore in 1993~\cite{moore1993}, 
while Chenciner and Montgomery proved the existence of this class of orbits a few years later~\cite{chenciner2000a}. 
In the early 2000s, many new choreographies for various $n$ and with unit masses were found by Sim\'{o} using a combination of different numerical methods~\cite{simo2001}. 
In 2008, Chen proved the existence of infinitely many choreographies of the three-body problem with a certain topological type for various choices of masses~\cite{chen2008}.

The situation is different on the sphere~\footnote{Throughout this paper, the term ``sphere'' will always refer to the 2-sphere.}. 
Although there is a growing interest in the $n$-body problem on the sphere (and in other spaces of constant curvature) in a cotangent potential~\cite{borisov2004,carinena2005,diacu2012c,diacu2013b,diacu2011,diacu2012a,diacu2012b,perez2012}, the only non-circular choreographies with unit masses found so far are for the two-body problem~\cite{borisov2004}~\footnote{Diacu proved the existence of choreographies with unit masses for the two- and three-body problems on the 3-sphere in~\cite{diacu2013c}.}.

We present in this paper an algorithm to compute planar and spherical choreographies with unit masses to high accuracy.
We have found many new choreographies on the sphere of radius $R$ in a cotangent potential for various $n\geq2$. 
These are curved versions of the planar choreographies found by Sim\'{o} and, as $R\rightarrow\infty$, they converge to the planar ones at a rate proportional to the curvature $1/R^2$.
 
\section{Planar choreographies}
 
Let $z_j(t)\in\C$, $0\leq j\leq n-1$, denote the positions of $n$ bodies with unit mass in the complex plane. 
The planar $n$-body problem describes the motion of these bodies under the action of Newton's law of gravitation, through the nonlinear coupled system of ODEs
\begin{equation}
\dsp z_j''(t) - \sum_{\underset{i\neq j}{i=0}}^{n-1} \frac{z_i(t) - z_j(t)}{\big\vert z_i(t) - z_j(t) \big\vert^3} = 0, \quad 0\leq j\leq n-1.
\label{newton}
\end{equation}

\noindent We are interested in periodic solutions of \reff{newton} in which the bodies share a single orbit and are uniformly spread along it, that is, solutions $z_j(t)$ such that
\begin{equation}
z_j(t) = q\Big(t + \frac{2\pi j}{n}	\Big), \quad 0\leq j\leq n-1,
\label{choreographies}
\end{equation}

\noindent for some $2\pi$-periodic function $q:[0,2\pi]\rightarrow\C$. Such solutions were named choreographies by Sim\'{o},
the $n$ bodies being ``seen to dance in a somewhat complicated way'' \cite{simo2001}. 
The period can be chosen equal to $2\pi$ because if $q(t)$ is a $T$-periodic solution of \reff{newton}, then $\lambda^{-2/3}q(\lambda t)$, $\lambda=T/(2\pi)$, is a $2\pi$-periodic one.
It has been well known since Poincar\'{e} \cite{poincare1892,poincare1896} that the \textit{principle of least action}, first introduced by Maupertuis in 1744 \cite{maupertuis1744}, can
be used to characterize periodic solutions of \reff{newton}:
choreographies \reff{choreographies} are minima of the \textit{action functional}, or simply \textit{action}, defined as the integral over one period of the 
kinetic minus the potential energy,
\begin{equation}
A = \int_0^{2\pi} \big(K(t) - U(t)\big)\,dt,
\label{action}
\end{equation}

\noindent with kinetic energy
\begin{equation}
\dsp K(t) = \frac{1}{2}\sum_{j=0}^{n-1} \big\vert z_j'(t) \big\vert^2 = \frac{1}{2}\sum_{j=0}^{n-1} 
\Big\vert q'\Big(t + \frac{2\pi j}{n}\Big) \Big\vert^2
\label{kineticenergy}
\end{equation}

\noindent and potential energy
\begin{equation}
\dsp U(t) = -\sum_{j=0}^{n-1}\sum_{i=0}^{j-1} \big\vert z_i(t) - z_j(t) \big\vert^{-1} = -\sum_{j=0}^{n-1}\sum_{i=0}^{j-1}
\Big\vert q\Big(t + \frac{2\pi i}{n}\Big) - q\Big(t + \frac{2\pi j}{n}\Big) \Big\vert^{-1}.
\label{newtonpotential}
\end{equation}

\noindent Note that the action \reff{action} depends on $q(t)$ via $U(t)$ and on $q'(t)$ via $K(t)$. 
Since the integral of \reff{kineticenergy} does not depend on $j$ and the integral of \reff{newtonpotential} only depends on $i-j$,
the action functional can be rewritten
\begin{equation}
\dsp A = \frac{n}{2}\int_0^{2\pi} \big\vert q'(t) \big\vert^2 dt
+ \frac{n}{2}\sum_{j=1}^{n-1} \int_0^{2\pi} \Big\vert q(t) - q\Big(t + \frac{2\pi j}{n}\Big) \Big\vert^{-1}dt.
\label{action2}
\end{equation}

\noindent Planar choreographies correspond to functions $q(t)$ which minimize \reff{action2}. 

We are also interested in solutions of \reff{newton} in which the bodies share a single orbit $q(t)$ that is rotating 
with angular velocity $\omega$ relative to an inertial reference frame, i.e.,
\begin{equation}
z_j(t) = e^{i\omega t}q\Big(t + \frac{2\pi j}{n}\Big), \quad 0\leq j\leq n-1.
\label{relativechoreographies}
\end{equation}

\noindent Choreographies of the form \reff{relativechoreographies} are said to be \textit{relative}, as opposed to the \textit{absolute} choreographies \reff{choreographies}.
The action associated with relative planar choreographies is
\begin{equation}
\dsp A = \frac{n}{2}\int_0^{2\pi} \big\vert q'(t) + i\omega q(t)\big\vert^2 dt
+ \frac{n}{2}\sum_{j=1}^{n-1} \int_0^{2\pi} \Big\vert q(t) - q\Big(t + \frac{2\pi j}{n}\Big) \Big\vert^{-1}dt.
\label{action3}
\end{equation}

\noindent Note that \reff{action2} is the special case of \reff{action3} with $\omega=0$. 

\section{Computing planar choreographies}

Our method for computing planar choreographies is based on the minimization of the action \reff{action3} and uses two key ingredients:

{\em Ingredient 1. Trigonometric interpolation.} The function $q(t)$ is represented by its trigonometric interpolant in the $\exp(ikt)$ basis. 
The optimization variables are the real and imaginary parts of its Fourier coefficients. The action is computed with the exponentially 
accurate trapezoidal rule. 

{\em Ingredient 2. Closed-form expressions for the gradient and the Hessian.} Formulas for the gradient and the Hessian matrix of the action \reff{action3} with respect to the optimization 
variables are derived explicitly and used in the optimization algorithms.

The numerical optimization of the action is in two steps:

{\em Step 1. Quasi-Newton optimization methods.} Numerical optimization methods with the exact gradient and based on approximations of the Hessian are employed with a small number of optimization variables.
The accuracy of the solution at this stage is from one to five digits.
This step is computationally very cheap. 

{\em Step 2. Newton's method.} Once an approximation to a choreography has been computed via a quasi-Newton method, one can improve the accuracy to typically ten digits with a few steps of Newton's method with exact Hessian, and a larger number of optimization variables.
This step is computationally more expensive.

Let us start with a few words about the first ingredient.
The approach used by Sim\'{o} \cite{simo2001} is to decompose the function $q(t)$ into real and imaginary parts, and to represent each of them by a trigonometric interpolant in the $\sin(kt)$ and $\cos(kt)$ basis. 
In this paper, we use instead a trigonometric interpolant of the function $q(t)$ itself in the $\exp(ikt)$ basis. 
For an odd number $N$, let $\{t_j=2\pi j/N\}$, $0\leq j\leq N-1$, denote $N$ equispaced points in $[0,2\pi)$ and 
$\{q_j=q(t_j)\}$, $0\leq j\leq N-1$, the (complex) values of $q(t)$ at the $t_j$'s. The trigonometric interpolant $p_N(t)$ of $q(t)$ at these points is
defined by
\begin{equation}
\dsp p_N(t) = \sumodd c_k e^{ikt}, \quad t\in[0,2\pi],
\label{trigfun}
\end{equation}

\noindent with Fourier coefficients
\begin{equation}
\dsp c_k = \frac{1}{N} \sum_{j=0}^{N-1} q_je^{-ikt_j}, \quad\vert k\vert\leq\frac{N-1}{2}.
\label{trigcoeffs}
\end{equation}

\noindent The trigonometric interpolant problem goes back at least to the young Gauss's calculations of the orbit of the asteroid Ceres in 1801---it
 seems that planetary orbits and trigonometric interpolation share a long and on-going relationship.
Throughout this paper, the number of grid points $N$ will always be odd. 
All our results have analogues for $N$ even, but the formulas are different, and little would be gained by writing everything twice. 
If we replace $q(t)$ by its trigonometric interpolant \reff{trigfun}--\reff{trigcoeffs} with $c_k=u_k + i v_k$, the action \reff{action3} becomes a function of
the $2N$ real variables $\{u_k,v_k\}$, $\vert k\vert\leq(N-1)/2$. We are looking for solutions $q(t)$ without collisions. 
The integrands in \reff{action3} are therefore analytic and the trapezoidal rule converges exponentially \cite{trefethen2014}. 
We use Chebfun~v5.2.1~\cite{chebfun} to compute trigonometric interpolants. 
Chebfun is an open-source package, MATLAB-based, for computing with functions to 16-digit accuracy. 
Its recent extension to periodic functions \cite{wright2015} provides a very convenient framework for working with closed curves in the complex plane.

Let us now say more about the second ingredient. The exact gradient and exact Hessian are derived in Appendix~A.
The gradient can be computed in $O(nN^2)$ operations, while the computation of the Hessian requires $O(nN^3)$ operations. 

The numerical optimization of the action, in MATLAB~R2015b, is carried out in two steps. 
First, we apply a quasi-Newton method \cite[Chapter 6]{nocedal2006} using the exact gradient, and with a small number of Fourier coefficients, $N=55$ and $75$ in our experiments.
Quasi-Newton methods are based on the approximation of the Hessian matrix (or its inverse) using rank-one or -two updates specified by gradient evaluations. 
In MATLAB, the \texttt{fminunc} command implements various quasi-Newton methods and, among them, we choose the BFGS algorithm \cite{shanno1970}.
We take $O(N)$ iterations of the BFGS algorithm. 
The cost of this first step is thus $O(nN^3)$ since, at each iteration, BFGS computes the gradient in $O(nN^2)$ operations and
matrix-vector products in $O(N^2)$ operations.
Second, we perform a small number $s$ of iterations of an approximate Newton method with exact Hessian and $M$ Fourier coefficients, $M>N$. 
The starting point of the approximate Newton method is the output of the BFGS algorithm, padded with $M-N$ zeros.
An exact Newton method would have a $O(snM^3)$ cost, since, at each iteration, it requires
the computation of the exact Hessian ($O(nM^3)$ operations) and the solution a linear system ($O(M^3)$ operations).
To reduce this cost, we use an approximate Newton method. We compute the exact Hessian at the first iteration only, and compute its $LDL^T$ decomposition ($O(M^3)$ operations, a generalization of Cholesky
decomposition for symmetric matrices which are not positive definitive, typically half as expensive as $LU$ factorization). 
The first iteration has thus a $O(nM^3)$ cost.
The subsequent iterations of Newton's method do not recompute the Hessian but use this factorization instead. The solution of the linear system can then be computed in $O(M^2)$ operations at each iteration.
The total cost of this second step is then $O(nM^3)$, and the total cost of the optimization is $O(n(N^3 + M^3))$.

Let us add four comments about this optimization process. 
First, since the initial guess of Newton's method---the output of BFGS---is a good approximation of a choreography, 
the Hessian matrix does not vary significantly from one iteration to another. As a consequence, using the $LDL^T$ factorization of the Hessian of the first iterate at each iteration does not affect the convergence very much.
Second, at a minimum of the action, i.e., a choreography, the Hessian is positive definite, so we could in principle use the Cholesky decomposition instead of the $LDL^T$ decomposition. However, in practice, because the Hessian is computed at an approximation of a choreography, it often has some small negative eigenvalues. 
Third, we do not use Newton's method with exact Hessian from the beginning because it only converges for initial guesses close enough to the solution. 
Fourth, another option would be to use $M$ coefficients with the quasi-Newton method directly. However, we found in practice that BFGS with $M$ coefficients typically achieves an accuracy of $6$ digits at most, while Newton's method achieves an accuracy of 10 digits.

For both steps of the optimization, the accuracy is defined as the the $2$-norm of the residual of \reff{newton} divided by the $2$-norm of the solution (relative $2$-norm).  
The residual is computed in Chebfun with the \texttt{chebop} class~\cite{driscoll2008}, the Chebfun automatic solver of differential equations.
We also check that the $2$-norm of the gradient divided by the $2$-norm of the gradient of the initial guess (relative $2$-norm) is close to zero, and that the Fourier coefficients of the solution decay to sufficiently small values.

\begin{figure}
\centering
\includegraphics [scale=.5]{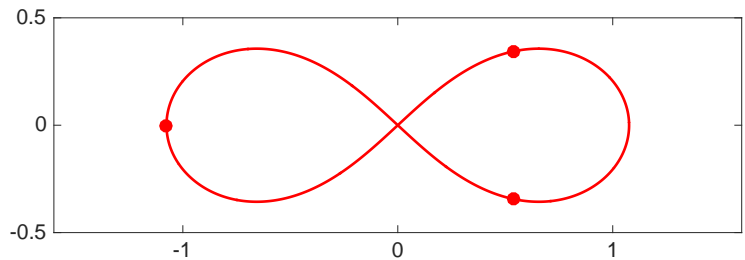}
\caption{\textit{The famous figure-eight solution of the three-body problem, obtained by running the code of Figure~$\ref{code}$. The dots show the bodies at
time $t=0$. The action of the resulting choreography, $A=24.371926476242812$, agrees with the $8$ digits given in \cite{chenciner2000a}.}}
\label{figure1}
\end{figure}

\begin{figure}
\centering
\includegraphics [scale=.37]{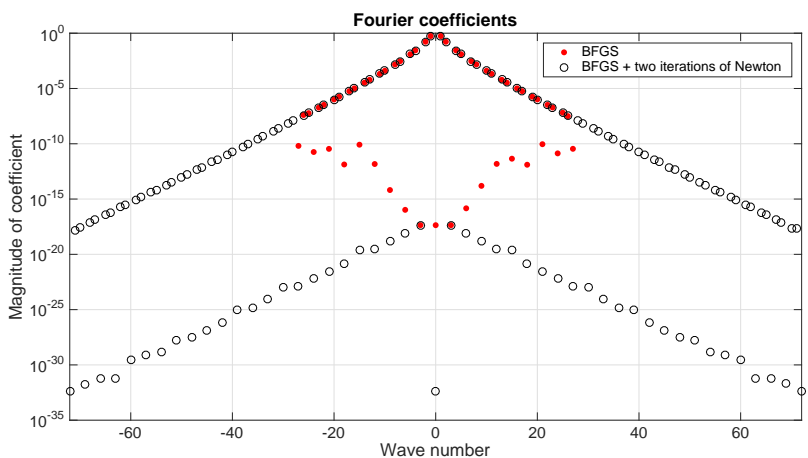}
\caption{\textit{Absolute values of the Fourier coefficients of the figure-eight of Figure~$\ref{figure1}$, obtained by BFGS (red dots) and BFGS followed by two iterations of Newton's method (black circles).}}
\label{figure1bis}
\end{figure}

The famous figure-eight, with action $A\approx24.371926$ \cite{chenciner2000a}, is plotted in Figure~\ref{figure1}. It is obtained by running the code of Figure~\ref{code}.
The code uses the \texttt{actiongradeval} and \texttt{gradhesseval} functions, which compute the action and the gradient, and the gradient and the Hessian; the codes are available online at the first author's GitHub web-page
 (\url{http://github.com/Hadrien-Montanelli}).
Table~\ref{table:figure1} shows some numbers pertaining to the computation of the figure-eight, including the relative $2$-norms of the solution and the gradient, and the amplitude of the smallest (numerically nonzero) Fourier coefficient.
After 51 iterations of the BFGS algorithm, the solution is accurate to five digits, and two iterations of Newton's method gives six extra correct digits. The solution found by BFGS and BFGS plus Newton look the same to the eye; if they were plotted on the same graph, they would be perfectly superimposed. The difference is visible in coefficient space. We plot the Fourier coefficients of the two solutions in Figure~\ref{figure1bis}. 
Since choreographies are analytic functions, the Fourier coefficients decay geometrically \cite[Theorem 4.1]{wright2015}.
The solution obtained by the BFGS algorithm only uses 55 coefficients, that is, wavenumbers $\vert k\vert\leq27$, and the Fourier coefficients decay to about $10^{-8}$. 
The solution obtained by Newton's method uses 145 Fourier coefficients, i.e., $\vert k\vert\leq72$, which decay to machine precision. 

\begin{table}
\vspace{.4cm}
\begin{center}
\begin{tabular}{l|cc}
& BFGS & Newton\\
\hline
Action & 24.371926476245442 & 24.371926476242809\\
Number of coefficients & 55 & 145\\
Computer time (s) & 0.19 & 0.27\\
Number of iterations & 51 & 2\\
Relative $2$-norm of the gradient & 2.86e-07 & 3.90e-16\\
Smallest coefficient & 3.17e-08 & 2.09e-18\\
Relative $2$-norm of the residual & 2.06e-05 & 2.24e-11\\
\end{tabular}
\end{center}
\vspace{.2cm}
\caption{\textit{Computation of the figure-eight choreography of Figure~$\ref{figure1}$.}}
\label{table:figure1}
\end{table}

\begin{figure}
\begin{verbatim}
% Initial guess:
  n = 3; N = 55; M = 145; 
  q0 = chebfun(@(t)cos(t)+1i*sin(2*t),[0 2*pi],N,'trig'); 
  c0 = trigcoeffs(q0);

% BFGS algorithm:
  options = optimoptions('fminunc');
  options.GradObj = 'on'; 
  options.Algorithm = 'quasi-newton'; 
  options.HessUpdate = 'bfgs';
  c = fminunc(@(x)actiongradeval(x,n),[real(c0);imag(c0)],options); 

% Two iterations of an approximate Newton method:
  c = [zeros((M-N)/2,1);c(1:N);zeros(M-N,1);c(N+1:end);zeros((M-N)/2,1)];
  mid = 1 + floor(M/2);
  [G, H] = gradhesseval(c,n); [L, D] = ldl(H);
  for k = 1:2
    s = L'\(D\(L\(-G)));
    cnew = c + [s(1:mid-1);0;s(mid:M+mid-2);0;s(M+mid-1:end)]; c = cnew;
    G = gradhesseval(c,n);
  end
  
% Reconstruct solution:
  c = c(1:M) + 1i*c(M+1:2*M);
  q = chebfun(c,[0 2*pi],'coeffs','trig');
\end{verbatim}
\caption{\textit{MATLAB code to compute the figure-eight of Figure~$\ref{figure1}$.
This code gives correct results to $11$ digits of accuracy in less than half a second on a $2.7$\,{\normalfont GHz} Intel i$7$ machine.}}
\label{code}
\end{figure}

\begin{figure}
\hspace{-1.3cm}
\includegraphics[scale=.4]{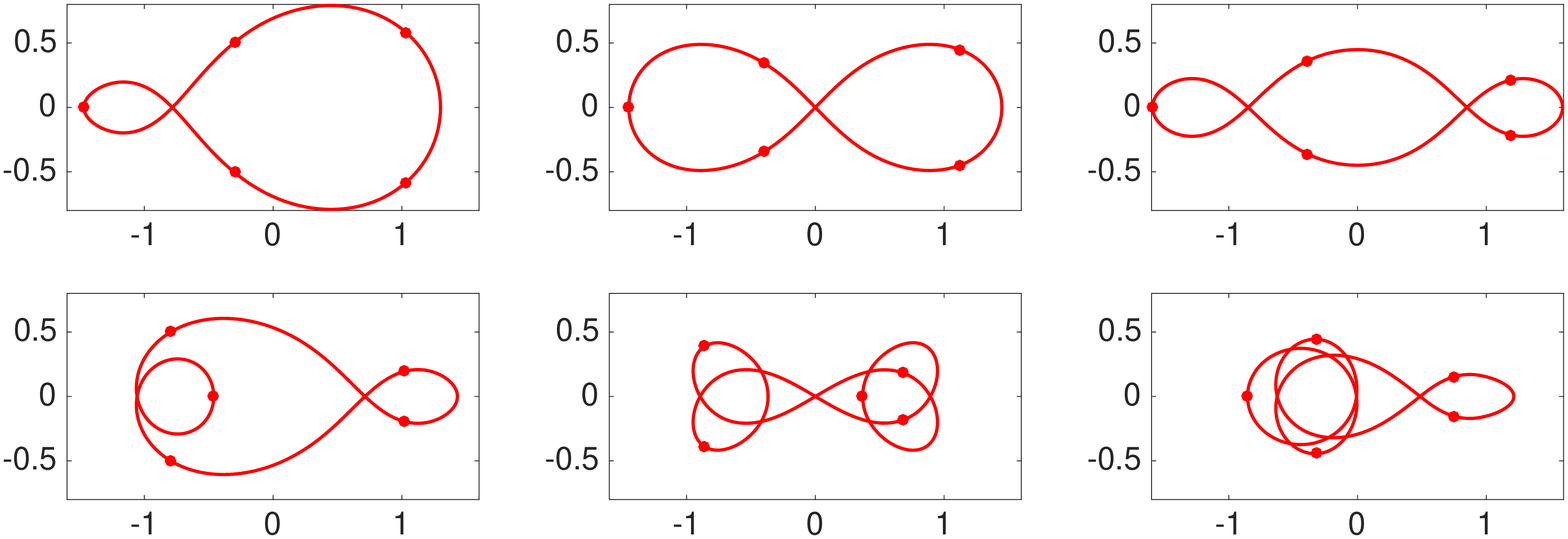}
\caption{\textit{Some absolute planar choreographies of the five-body problem. These correspond to choreographies~$1$ (top-left corner), $2$ (top center), $4$ (top-right corner), 
$7$ (bottom-left corner), $16$ (bottom center) and $18$ (bottom-right corner) of \cite[Figure 2]{simo2001}. 
The dots show the bodies at $t=0$.}}
\label{figure2}
\end{figure}

All the planar absolute choreographies of the five-body problem found by Sim\'{o} \cite[Figure 2]{simo2001} can be computed with this algorithm. 
We plot six of them in Figure~\ref{figure2}.
Tables~\ref{table1:figure2} and \ref{table2:figure2} show some numbers pertaining to their computation.
The BFGS algorithm leads to results accurate to a few digits, and two to six iterations of Newton's method lead to about ten digits of accuracy.

\begin{table}
\begin{center}
\begin{tabular}{l|cccccc}
& 1 & 2 & 4 & 7 & 16 & 18\\
\hline
Action &  68.8516 &  71.3312 & 77.1588 & 88.4397 & 109.6361 & 119.3168\\
Number of coefficients & 75 & 75 & 75 & 75 & 75 & 75\\
Computer time (s) & 0.89 & 0.63 & 0.59 & 0.62 & 1.73 & 0.83\\
Number of iterations & 98 & 68 & 69 & 136 & 373 & 166\\
Relative $2$-norm of the gradient & 3.27e-08 & 6.65e-08 & 3.36e-07 & 4.85e-09 & 3.21e-10 & 2.89e-09\\
Smallest coefficient                      & 1.28e-05 & 6.13e-09 & 5.89e-06 & 8.12e-06 & 3.61e-05 & 6.52e-05\\
Relative $2$-norm of the residual  & 3.51e-02 & 1.07e-05 & 6.65e-03 & 2.23e-02 & 2.75e-01 & 3.55e-01\\
\end{tabular}
\end{center}
\vspace{.2cm}
\caption{\textit{Computation of the absolute planar choreographies of Figure~$\ref{figure2}$ with the BFGS algorithm.}}
\label{table1:figure2}
\end{table}

\begin{table}[H]
\begin{center}
\begin{tabular}{l|cccccc}
& 1 & 2 & 4 & 7 & 16 & 18\\
\hline
Action &  68.8516 &  71.3312 & 77.1588 & 88.4397 & 109.6361 & 119.3184\\
Number of coefficients & 335 & 145 & 245 & 285 & 445 & 455\\
Computer time (s) & 3.49 & 0.48 & 1.68 & 2.33 & 7.26 & 7.85\\
Number of iterations & 4 & 2 & 3 & 4 & 5 & 6\\
Relative $2$-norm of the gradient & 6.93e-16 & 1.90e-16 & 1.71e-13 & 1.53e-16 & 3.74e-14 & 3.79e-14\\
Smallest coefficient                      & 3.43e-15 & 8.81e-16 & 8.00e-17 & 2.05e-14 & 1.26e-16 & 1.43e-14\\
Relative $2$-norm of the residual  & 8.29e-10 & 5.66e-11 & 7.27e-10 & 8.31e-10 & 1.51e-09 & 2.62e-09\\
\end{tabular}
\end{center}
\vspace{.2cm}
\caption{\textit{Computation of the absolute planar choreographies of Figure~$\ref{figure2}$ with Newton's method using the outputs of BFGS as initial guesses.}}
\label{table2:figure2}
\end{table}

\begin{figure}
\hspace{-.5cm}
\includegraphics[scale=.35]{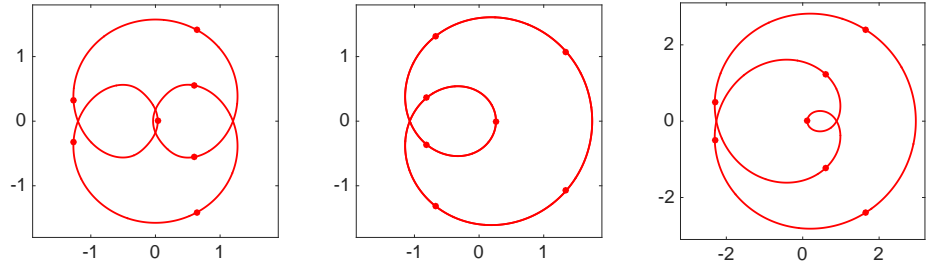}
\caption{\textit{Relative planar choreographies of the seven-body problem with angular velocity $2.8$ (left), $-2.9$ (center) and $2.31$ (right).
The dots show the bodies at time $t=0$. They can be computed to ten digits of accuracy with about four hundred Fourier coefficients.}}
\label{figure3}
\end{figure}

The same method can be used to compute relative choreographies to high accuracy. We plot three relative planar choreographies of the seven-body problem in Figure~\ref{figure3}.

An interactive tool to compute choreographies with MATLAB and Chebfun is available at the web-page previously given.
The code, \texttt{choreo}, finds choreographies starting with hand-drawn initial guesses. It is easy to use, fast and enjoyable---the reader is highly encouraged to try it!

Let us conclude this section with a few words about the number of choreographies for a given $n$. 
This number is not known, but there is an interesting result, due to Sim\'{o} \cite[Proposition 5.1]{simo2001}, about the (smaller) number of choreographies that consist of
a concatenation of ``bubbles,'' such as choreographies 1, 2 and 4 of Figure~\ref{figure2}. For $n\geq3$, there are $2^{n-3} + 2^{\lfloor (n-3)/2\rfloor}$ such choreographies.

\section{Spherical choreographies}

Let $X_j(t)\in\R^3$, $0\leq j\leq n-1$, denote the Cartesian coordinates of $n$ bodies with unit mass on the sphere 
$\Stwo_R=\{X\in\R^3, \Vert X\Vert=R\}$, where $\Vert\cdot\Vert$ is the Euclidean norm in $\R^3$.
The $n$-body problem on the sphere in a cotangent potential describes the motion of these bodies via the $n$ coupled nonlinear ODEs 
\begin{equation}
\dsp X_j''(t) - \sum_{\underset{i\neq j}{i=0}}^{n-1}\frac{R^3X_i(t) - R(X_i(t)\cdot X_j(t))X_j(t)}{\big[R^4-(X_i(t)\cdot X_j(t))^2\big]^{3/2}}
+ R^{-2}\big\Vert X_j'(t) \big\Vert^2 X_j(t) = 0, \; 0\leq j\leq n-1.
\label{sphericalnewton}
\end{equation}

\noindent See \cite{diacu2012a} for details about the derivation of these equations.
Note that the potential associated with \reff{sphericalnewton} is no longer the Newtonian potential \reff{newtonpotential}. 
It is a cotangent potential, a generalization of the Newtonian potential on the sphere, and dates back to the 1820's with the work
of Bolyai and Lobachevsky. The reader can find a detailed history of the problem in Diacu's 2012 book on relative equilibria \cite{diacu2012c}.

We are looking for periodic solutions of \reff{sphericalnewton} moving along the same orbit, i.e., solutions $X_j(t)$ such that
\begin{equation}
X_j(t) = Q\Big(t + \frac{2\pi j}{n}	\Big), \quad 0\leq j\leq n-1,
\label{sphericalchoreographies}
\end{equation}

\noindent for some $2\pi$-periodic function $Q:[0,2\pi]\rightarrow\Stwo_R\subset\R^3$. 
Again, the period can be chosen equal to $2\pi$ because if $Q(t)$ is a $T$-periodic solution of \reff{sphericalnewton} on the sphere of radius $R$, 
then $\lambda^{-2/3}Q(\lambda t)$, $\lambda=T/(2\pi)$, is a $2\pi$-periodic one on the sphere of radius $\lambda^{-2/3}R$.
We call these solutions \textit{spherical choreographies}.
They are minima of the action associated with \reff{sphericalnewton}, defined again as the integral over one period of the kinetic minus the potential energy,
with kinetic energy
\begin{equation}
\dsp K(t) = \frac{1}{2}\sum_{j=0}^{n-1} \big\Vert X_j'(t) \big\Vert^2 = \frac{1}{2}\sum_{j=0}^{n-1} 
\Big\Vert Q'\Big(t + \frac{2\pi j}{n}\Big) \Big\Vert^2
\end{equation}

\noindent and potential energy
\begin{equation}
\dsp U(t) = -\frac{1}{R}\sum_{j=0}^{n-1}\sum_{i=0}^{j-1} \cot \frac{\hat{d}(X_i(t),X_j(t))}{R},
\label{cotangentpotential}
\end{equation}

\noindent where 
\begin{equation}
\hat{d}(X_i(t),X_j(t)) = R\arccos\frac{X_i(t)\cdot X_j(t)}{R^2}
\label{greatcircledistance}
\end{equation}

\noindent is the great-circle distance between $X_i(t)$ and $X_j(t)$ on $\Stwo_R$. The potential \reff{cotangentpotential} is the cotangent of the (rescaled) distance on the sphere.
Using the trigonometric identity $\cot(\arccos(x))=x/\sqrt{1-x^2}$, the potential energy can be rewritten
\begin{equation}
\dsp U(t) = -\frac{1}{R}\sum_{j=0}^{n-1}\sum_{i=0}^{j-1}\frac{X_i(t)\cdot X_j(t)}{\sqrt{R^4-(X_i(t)\cdot X_j(t))^2}}.
\end{equation}

\noindent The action is then given by
\begin{equation}
\dsp A = \frac{n}{2}\int_0^{2\pi} \big\Vert Q'(t)\big\Vert^2 dt
+ \frac{n}{2R}\sum_{j=1}^{n-1}\int_0^{2\pi} \frac{Q(t)\cdot Q\big(t + \frac{2\pi j}{n}\big)}{\sqrt{R^4-\big(Q(t)\cdot Q\big(t + \frac{2\pi j}{n}\big)\big)^2}}\,dt.
\label{sphericalaction}
\end{equation}

\noindent Spherical choreographies correspond to functions $Q(t)$ which minimize \reff{sphericalaction}. 
Note that since the cotangent potential \reff{cotangentpotential} is singular not only when the distance between two bodies is zero 
but also for antipodal configurations, we are looking for solutions that stay in a single hemisphere. 
See \cite{diacu2012b} for more details about the singularities of the $n$-body problem in a cotangent potential.

As in the plane, we are also interested in solutions of \reff{sphericalnewton} in which the bodies share a single orbit $Q(t)$ that is 
rotating with angular velocity $\omega$ along the $z$-axis relative to an inertial reference frame, i.e., 
\begin{equation}
X_j(t) = \begin{bmatrix}
\cos(\omega t) & -\sin(\omega t) & 0\\
\sin(\omega t) & \cos(\omega t) & 0\\
0 & 0 & 1
\end{bmatrix}
Q\Big(t + \frac{2\pi j}{n}	\Big), \quad 0\leq j\leq n-1.
\label{relativesphericalchoreographies}
\end{equation}

\noindent Let $R_\omega(t)$ denote the rotation matrix in \reff{relativesphericalchoreographies}. 
The action associated with relative spherical choreographies is
\begin{equation}
\dsp A = \frac{n}{2}\int_0^{2\pi} \big\Vert R_\omega(t)Q'(t) + R_\omega'(t)Q(t)\big\Vert^2 dt
+ \frac{n}{2R}\sum_{j=1}^{n-1}\int_0^{2\pi} \frac{Q(t)\cdot Q\big(t + \frac{2\pi j}{n}\big)}{\sqrt{R^4-\big(Q(t)\cdot Q\big(t + \frac{2\pi j}{n}\big)\big)^2}}\,dt.
\label{sphericalaction2}
\end{equation}

\noindent As in the plane, \reff{sphericalaction} is the special case of \reff{sphericalaction2} with $\omega=0$.

\section{Computing spherical choreographies}

Our method for computing spherical choreographies is based on \textit{stereographic projection} and on the algorithm described in Section~3. 
Points $X=(x_1,x_2,x_3)^T$ on the sphere $\Stwo_R$ are mapped to points $z=P_R(X)$ in the plane $\C$ via
\begin{equation}
z = P_R(X) = \frac{R x_1 + i R x_2}{R - x_3}.
\label{stereoproj}
\end{equation}

\noindent The inverse mapping is given by
\begin{equation}
X = P_R^{-1}(z) = \frac{1}{R^2+\vert z\vert^2}(2R^2\mrm{Re}(z),2R^2\mrm{Im}(z),-R^3 + R\vert z\vert^2)^T.
\label{stereoprojinv}
\end{equation}

\noindent The Euclidean distance $d(X,Y)=\Vert X-Y\Vert$ between two points on the sphere is transformed into the distance $d(z,\xi)$ between their projections
$z=P_R(X)$ and $\xi=P_R(Y)$ defined by
\begin{equation}
d(z,\xi) = \frac{2R^2\vert z - \xi \vert}{\sqrt{(R^2+\vert z\vert^2)(R^2+\vert \xi\vert^2)}},
\label{distancebis}
\end{equation}

\noindent and the great-circle distance \reff{greatcircledistance} into
\begin{equation}
\hat{d}(z,\xi) = 2R\arcsin\frac{d(z,\xi)}{2R}.
\label{greatcircledistancebis}
\end{equation}

\noindent The complex plane endowed with the distance \reff{greatcircledistancebis} is called the \textit{spherical plane}. 
Let $q(t)=P_R(Q(t))$ denote the projection of the curve $Q(t)$ onto $\C$, and 
\begin{equation}
z_j(t)=P_R(X_j(t))=P_R\Big(Q\Big(t+\frac{2\pi j}{n}\Big)\Big)=q\Big(t+\frac{2\pi j}{n}\Big), \quad 0\leq j\leq n-1,
\end{equation}

\noindent the projections of the $n$ bodies $X_j(t)$. The action \reff{sphericalaction2} can be then reformulated as
\begin{equation}
\dsp A =  \dsp\frac{n}{2}\int_0^{2\pi}\bigg(\frac{2R^2\vert q'(t) + i\omega q(t)\vert}{R^2+\vert q(t)\vert^2}\bigg)^2dt
+ \frac{n}{2R}\sum_{j=1}^{n-1}\int_0^{2\pi}\frac{2R^2 - D_j(t)^2}{D_j(t)\sqrt{4R^2-D_j(t)^2}}\,dt,
\label{sphericalaction2bis}
\end{equation}
 
\noindent with $D_j(t) = d\big(q(t),q\big(t+\frac{2\pi j}{n}\big)\big)$.
P\'{e}rez-Chavela and Reyes-Victoria \cite[Theorem 2.3]{perez2012} showed the equivalence of the formulations \reff{sphericalaction2} and \reff{sphericalaction2bis} with $\omega=0$.
 
Once the problem is reformulated in the spherical plane, we apply the two key ingredients and the two steps described in Section~3. 
The function $q(t)$ is approximated by its trigonometric interpolant \reff{trigfun}--\reff{trigcoeffs} at $N$ points, the action \reff{sphericalaction2bis} becomes a function of the real and imaginary parts of the Fourier coefficients, and is computed with the exponentially accurate trapezoidal rule.
Formulas for the gradient and the Hessian matrix of the action \reff{sphericalaction2bis} are derived in Appendix~B and used in the optimization algorithms. Again, the computation of the gradient costs $O(nN^2)$ while the computation of the Hessian requires $O(nN^3)$ operations.
For the optimization, we use the same strategy: BFGS algorithm with exact gradient and a small number of variables, followed by a few steps of an approximate Newton method with exact Hessian and a larger number of variables. As in the plane, at convergence, we check that the norm of the gradient of the action is close to zero, the Fourier coefficients of the solution decay to sufficiently small values, and the solution satisfies equation \reff{sphericalnewton} projected into the plane. 
The latter was first given by P\'{e}rez-Chavela and Reyes-Victoria in 2012 \cite[Lemma 2.1]{perez2012}, and can be written as
\begin{equation}
z_j''(t) = \frac{2\bar{z}_j(t)z_j'^2(t)}{R^2+\vert z_j(t)\vert^2} + \frac{4R}{\lambda_j(t)}\sum_{\underset{i\neq j}{i=0}}^{n-1}\frac{P_{j,i}(t)}{\Theta_{j,i}(t)^{3/2}}, \quad 0\leq j\leq n-1,
\label{sphericalnewtonbis}
\end{equation}

\noindent where $\lambda_j(t)=4R^4/(R^2+\vert z_j(t)\vert^2)^2$ is the conformal factor that appears in the kinetic part of \reff{sphericalaction2bis}, while $P_{j,i}(t)$ and $\Theta_{j,i}(t)$ are defined by
\begin{equation}
P_{j,i}(t) = \big[R^2+\vert z_j(t)\vert^2\big]\big[R^2+\vert z_i(t)\vert^2\big]^2\big[R^2+\bar{z}_i(t)z_j(t)\big]\big[z_i(t)-z_j(t)\big],
\end{equation}

\noindent and 
\begin{equation}
\begin{array}{ll}
\Theta_{j,i}(t) & = - \big[2R^2z_j(t)\bar{z}_i(t) + 2R^2z_i(t)\bar{z}_j(t) + (\vert z_j(t)\vert^2-R^2)(\vert z_i(t)\vert^2-R^2)\big]^2\\
& + \big[R^2+\vert z_j(t)\vert^2\big]^2\big[R^2+\vert z_i(t)\vert^2\big]^2.
\end{array}
\end{equation}

\noindent Again, the residual of equation \reff{sphericalnewtonbis} can be computed in Chebfun with \texttt{chebop}.

\begin{figure}
\centering
\includegraphics [scale=.42]{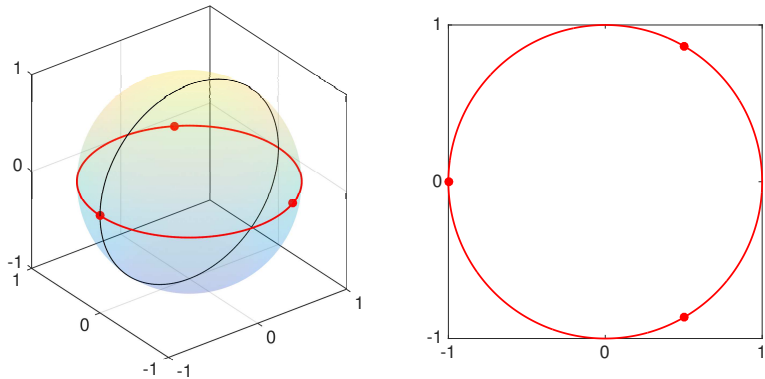}
\caption{\textit{A circular choreography of the spherical three-body problem on the sphere of radius $R=1$ (left) and its projection in the plane (right). 
All the circles of radius $r$, $0<r<R$, are spherical choreographies for any $R\neq0$ and any $n\geq2$.
Circles of radius $r=R$, i.e. equators, such as the one above, are spherical choreographies for odd $n$ only;
for even $n$, it would lead to antipodal singularities.
The dots show the bodies at time $t=0$.}}
\label{figure4}
\end{figure}

As we mentioned in the introduction, the only non-circular spherical choreographies with unit masses found so far are for the two-body problem.
Diacu and its collaborators \cite{diacu2012a}, and P\'{e}rez-Chavela and Reyes-Victoria \cite{perez2012} also characterized the solutions of the spherical $n$-body problem, $n\geq 2$, in which the bodies move
along the same circle (such as Figure~\ref{figure4}), or along different ones---the \textit{relative equilibria}. 

\begin{figure}
\hspace{-1.75cm}
\includegraphics [scale=.55]{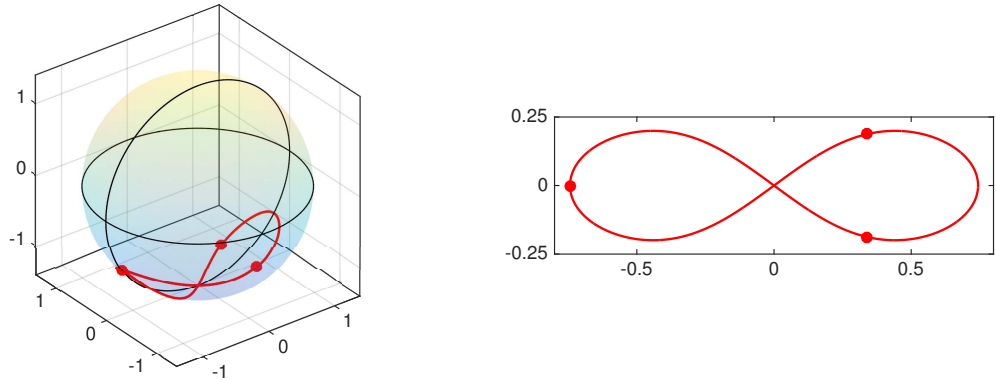}
\caption{\textit{Spherical figure-eight on the sphere of radius $1.4$ (left) and its projection in the plane (right).
The dots show the bodies at time $t=0$.}}
\label{figure5}
\end{figure}

\begin{figure}
\centering
\includegraphics [scale=.35]{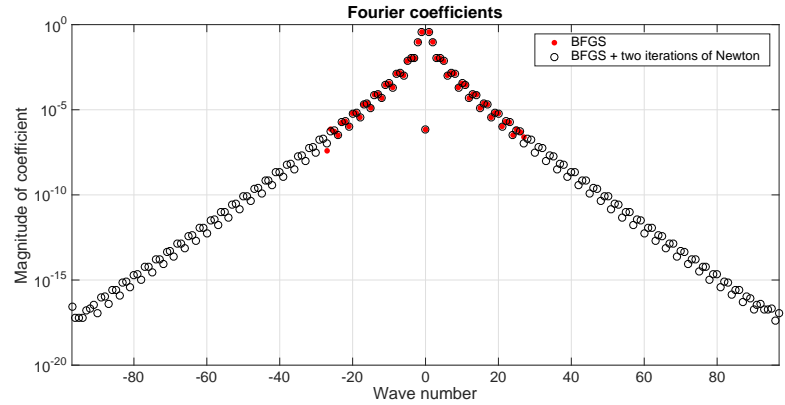}
\caption{\textit{Absolute values of the Fourier coefficients of the spherical figure-eight of Figure~$\ref{figure5}$, obtained by BFGS (red dots) and BFGS followed by one step of Newton's method (black circles).}}
\label{figure5bis}
\end{figure}

We present now new non-circular spherical choreographies.
The first one is the \textit{spherical figure-eight}, a solution of the three-body problem on the sphere of radius $R=1.4$, shown in Figure~\ref{figure5}. 
Table \ref{table:figure5} shows some numbers pertaining to its computation. Combining the BFGS algorithm with Newton's method leads to thirteen digits of accuracy.
We plot the geometrically decaying Fourier coefficients of the outputs of BFGS and Newton's method in Figure~\ref{figure5bis}. 
After BFGS, they decay to $10^{-7}$, and after two iterations of Newton's method, they decay to machine precision.
Numerically, we found that the ($2\pi$-periodic) spherical figure-eight exists on spheres of radius $R\geq 1.32$. 
Below this value, it cannot fit in a single hemisphere and would therefore lead to antipodal singularities 
\footnote{The solution of Figure~\ref{figure5} exists for radii $R<1.32$ but with a shorter period $T$, with $R^3/T^2=(1.4)^3/(2\pi)^2$.
This is a consequence of the scaling invariance described below equation \reff{sphericalchoreographies}.}.

\begin{table}
\vspace{.4cm}
\begin{center}
\begin{tabular}{l|cc}
& BFGS & Newton\\
\hline
Action & 18.948304138530286 & 18.948304135957898\\
Number of coefficients & 55 & 195\\
Computer time (s) & 0.51 & 3.19\\
Number of iterations & 72 & 2\\
Relative $2$-norm of the gradient & 5.18e-07 & 1.10e-13\\
Smallest coefficient & 2.48e-07 & 1.06e-17\\
Relative $2$-norm of the residual & 4.08e-04 & 8.82e-13\\
\end{tabular}
\end{center}
\vspace{.2cm}
\caption{\textit{Computation of the spherical figure-eight choreography of Figure~$\ref{figure5}$.}}
\label{table:figure5}
\end{table}

\begin{figure}
\hspace{-1.7cm}
\includegraphics [scale=.5]{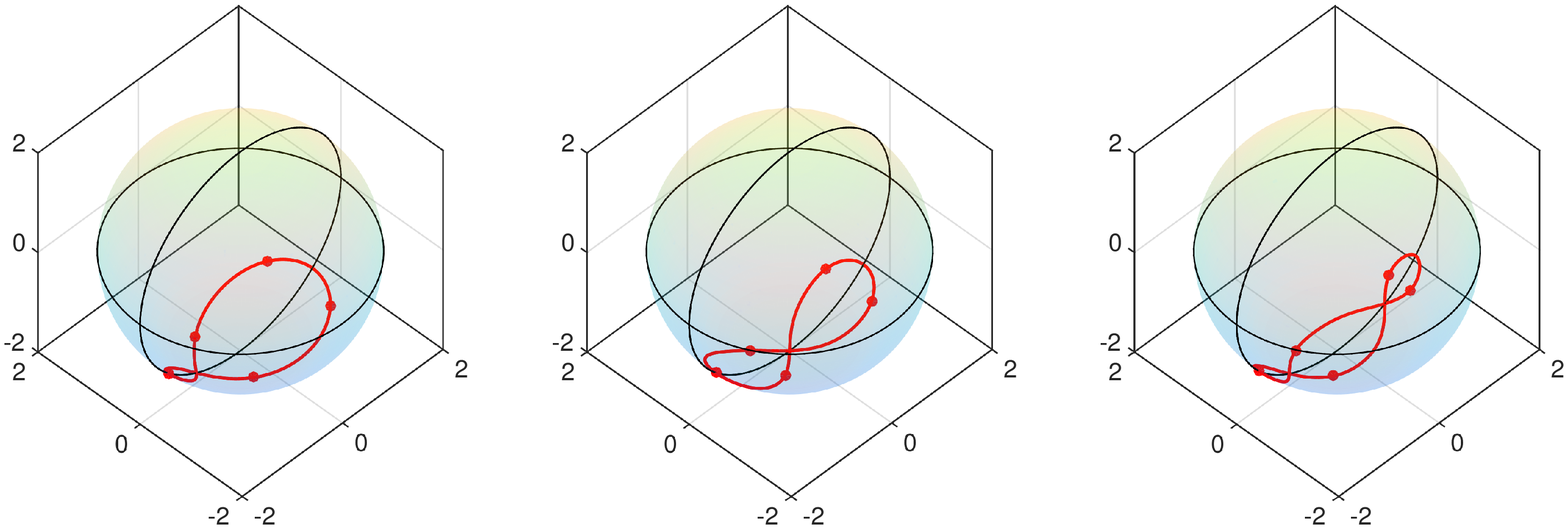}
\caption{\textit{Spherical choreographies of the five-body problem on the sphere of radius $2$, analogous to the planar choreographies of Figure~$\ref{figure2}$, with action 
$58.6829$ (left), $61.8035$ (center), and $67.4127$ (right).}}
\label{figure6}
\end{figure}

Many new spherical choreographies can be found with our algorithm. We show in Figure~\ref{figure6} three spherical choreographies
of the five-body problem on the sphere of radius 2. These are curved versions of the choreographies of Figure~\ref{figure2}. Table~\ref{table:figure6} shows some numbers pertaining to their computations. 
We get two to five digits of accuracy with the BFGS algorithm, and applying Newton's method with the outputs of BFGS as initial guesses leads to nine to thirteen digits of accuracy.

\begin{table}
\begin{center}
\begin{tabular}{l|cc|cc|cc}
& BFGS & Newton & BFGS & Newton & BFGS & Newton\\
\hline
Action & 58.6831 & 58.6829 & 61.8035 & 61.8035 & 67.4127 &  67.4127\\
Number of coefficients & 75 &  375 & 75 & 205 & 75 & 245\\
Computer time (s) & 2.35 & 33.13 & 1.05 & 6.08 & 2.14 & 9.44\\
Number of iterations & 142 & 6 & 65 & 3 & 115 & 4\\
Relative $2$-norm of the gradient & 3.58e-05 & 9.81e-11 & 3.50e-07 & 7.61e-14 & 1.98e-07 & 5.72e-14\\
Smallest coefficient                      & 1.88e-05 & 2.87e-14 & 4.52e-08 & 1.16e-17 & 6.37e-06 & 2.49e-15\\
Relative $2$-norm of the residual  & 4.26e-02 & 1.68e-09 & 5.62e-05 & 7.32e-13 & 7.05e-03 & 9.32e-09\\
\end{tabular}
\end{center}
\vspace{.2cm}
\caption{\textit{Computation of the spherical choreographies of Figure~$\ref{figure6}$.}}
\label{table:figure6}
\end{table}

\begin{figure}
\hspace{-1.7cm}
\includegraphics [scale=.5]{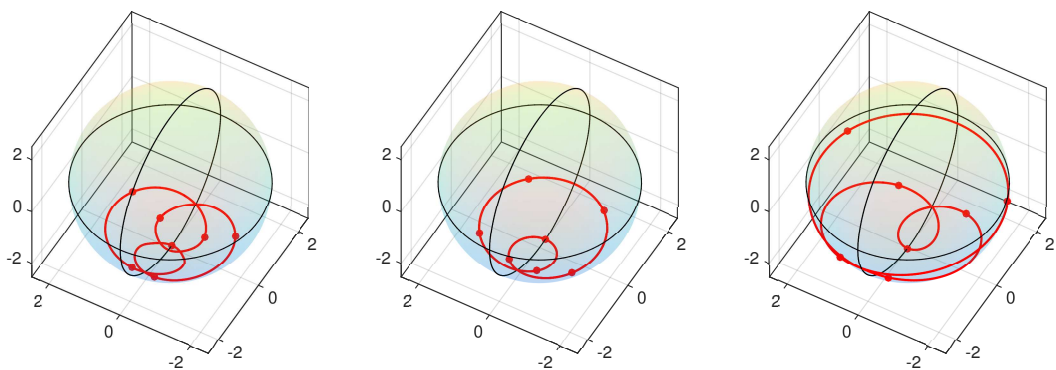}
\caption{\textit{Relative spherical choreographies of the seven-body problem on the sphere of radius $2.5$ with angular velocities $2.8$ (left), $-2.9$ (center) and $2.31$ (right),
analogous to the relative choreographies of Figure~$\ref{figure3}$. Again, they can be computed to ten digits of accuracy with about four hundred Fourier coefficients.}}
\label{figure7}
\end{figure}

Relative spherical choreographies can also be computed with this method. We plot three relative spherical choreographies of the seven-body problem on the sphere of radius 2.5 in Figure~\ref{figure7}.
These are curved versions of the relative choreographies of Figure~\ref{figure3}.

An interactive tool to compute spherical choreographies using hand-drawn initial guesses, \texttt{choreosphere}, is also available at the web-page previously given;
it uses the \texttt{actiongradevalsphere} and \texttt{gradhessevalsphere} functions, which compute the action and the gradient, and the gradient and the Hessian.

\section{Limit of infinitely large radius}

\begin{figure}
\hspace{1.4cm}
\includegraphics [scale=.35]{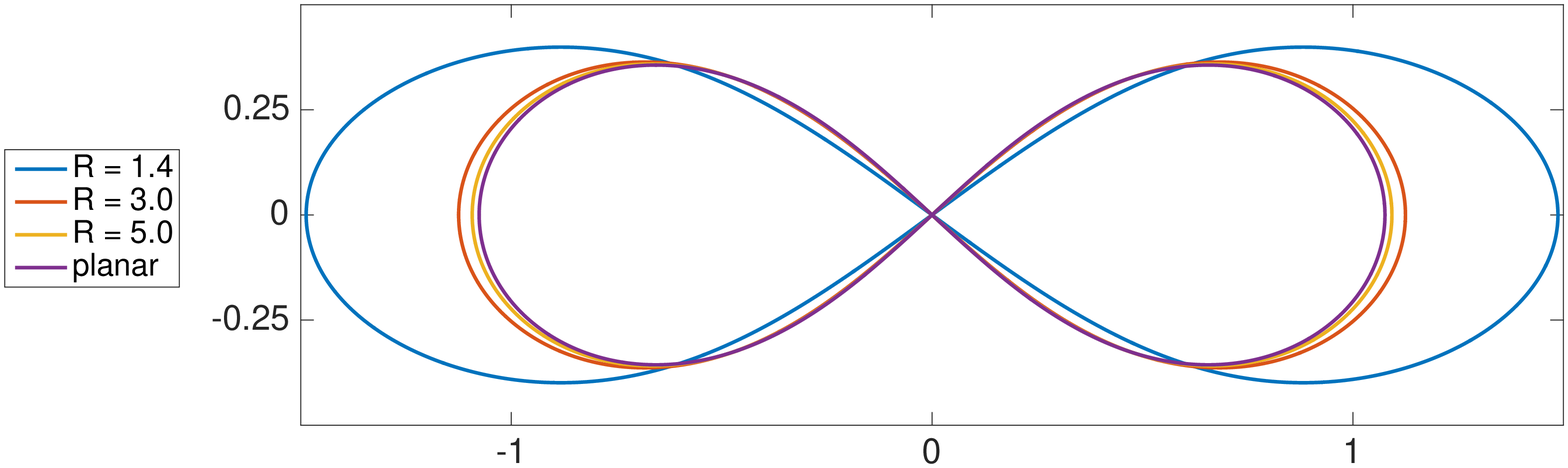}
\caption{\textit{Spherical figure-eight of Figure~$\ref{figure5}$ (multiplied by a factor $2$) for different values of $R$, together with its planar analogue of Figure~$\ref{figure1}$.
As $R$ increases, the spherical figure-eight converges to the planar one.}}
\label{figure5tri}
\end{figure}

\begin{figure}
\hspace{.05cm}
\includegraphics [scale=.37]{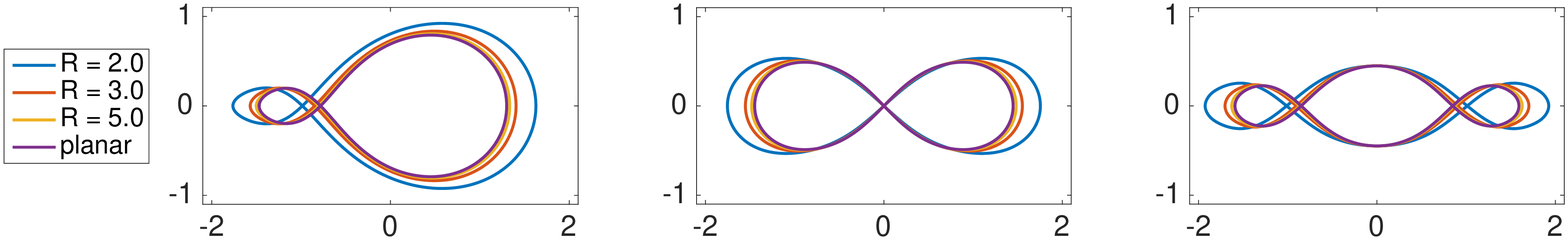}
\caption{\textit{Spherical choreographies of Figure~$\ref{figure6}$ (multiplied by a factor $2$) for different values of $R$, together with their planar analogues of Figure~$\ref{figure2}$.
As $R$ increases, the spherical choreographies converge to the planar ones.}}
\label{figure8}
\end{figure}

\begin{figure}
\hspace{.5cm}
\includegraphics [scale=.35]{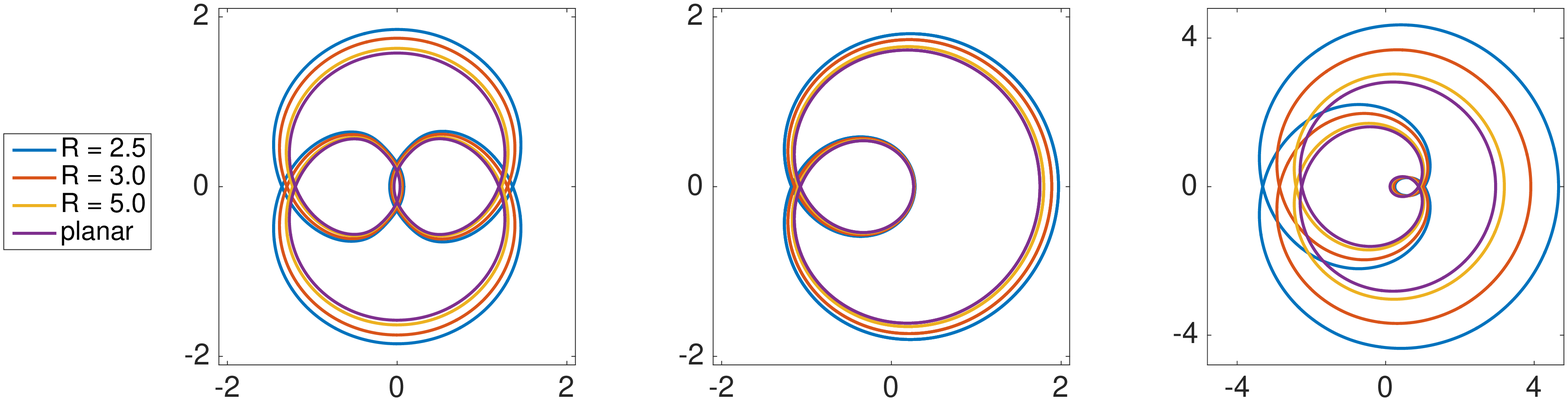}
\caption{\textit{Spherical choreographies of Figure~$\ref{figure7}$ (multiplied by a factor $2$) for different values of $R$, together with their planar analogues of Figure~$\ref{figure3}$.
As $R$ increases, the spherical choreographies converge to the planar ones.}}
\label{figure9}
\end{figure}

As its radius $R$ gets bigger, the sphere gets flatter, and in the limit $R\rightarrow\infty$, it converges to the complex plane. 
Equivalently, the spherical plane converges to the complex plane.
The distances \reff{distancebis} and \reff{greatcircledistancebis} converge to twice the absolute value, and the action on the sphere \reff{sphericalaction2bis} converges to four times
the action in the plane \reff{action3}, since it involves squares of distances.
We might then expect that twice the spherical choreographies converge to the planar choreographies as $R\rightarrow\infty$, and it is indeed the case 
\footnote{We are studying the convergence with a fixed period $2\pi$. Similarly, $T$-periodic spherical choreographies converge to $T$-periodic planar choreographies.}. 
In Figures~\ref{figure5tri}, \ref{figure8} and \ref{figure9}, we plot the spherical choreographies of Figures~\ref{figure5}, \ref{figure6} and \ref{figure7} (multiplied by a factor 2) for increasing values of $R$ and plot them together with their planar analogues.
Tables~\ref{table:figure5tri}, \ref{table:figure8} and \ref{table:figure9} report the $\infty$-norm of the difference between analogous spherical and planar choreographies as $R$ increases.
It is clear from the tables that spherical choreographies converge to their planar analogues at a rate proportional to the curvature $1/R^2$.

\begin{table}
\begin{center}
\begin{tabular}{l|cccccc}
& $R$ = 1.4 & 3 & 5 & 10 & 100 & 1000\\
\hline
& 4.11e-01 & 4.87e-02 & 1.65e-02 & 4.03e-03 & 4.00e-05 & 4.03e-07\\
\end{tabular}
\vspace{.2cm}
\caption{\textit{Convergence of the spherical figure-eight of Figure~$\ref{figure5}$ to the planar figure-eight of Figure~$\ref{figure1}$ as $R$ increases.}}
\label{table:figure5tri}
\end{center}
\end{table}

\begin{table}
\begin{center}
\begin{tabular}{l|cccccc}
& $R$ = 2 & 3 & 5 & 10 & 100 & 1000\\
\hline
Left & 3.27e-01 & 1.05e-01 & 3.37e-02 & 8.09e-03 & 7.98e-05 &  8.00e-07\\
Middle & 3.06e-01 & 1.04e-01 & 3.39e-02 & 8.17e-03 &  8.07e-05 & 8.04e-07\\
Right & 3.34e-01 & 1.12e-01 & 3.64e-02 & 8.75e-03 & 8.65e-05 & 8.65e-07\\
\end{tabular}
\vspace{.2cm}
\caption{\textit{Convergence of the spherical choreographies of Figure~$\ref{figure6}$ to the planar ones of Figure~$\ref{figure2}$ as $R$ increases.}}
\label{table:figure8}
\end{center}
\end{table}

\begin{table}
\begin{center}
\begin{tabular}{l|cccccc}
& $R$ = 2.5 & 3 & 5 & 10 & 100 & 1000\\
\hline
Left & 2.78e-01 & 1.74e-01 & 5.54e-02 & 1.32e-02 & 1.30e-04 &  1.30e-06\\
Middle & 2.19e-01 & 1.39e-01 & 4.50e-02 & 1.08e-02 &  1.07e-04 & 1.12e-06\\
Right & 1.69e+00 & 9.47e-01 & 2.35e-01 & 5.22e-02 & 5.04e-04 & 5.02e-06\\
\end{tabular}
\vspace{.2cm}
\caption{\textit{Convergence of the spherical choreographies of Figure~$\ref{figure7}$ to the planar ones of Figure~$\ref{figure3}$ as $R$ increases.}}
\label{table:figure9}
\end{center}
\end{table}

\section{Conclusions}

Choreographies are very special solutions of the $n$-body problem. They are not only periodic but also share a single orbit.
We have shown in this paper that choreographies exist on a sphere in a cotangent potential for various $n\geq2$. 
Curved versions of Sim\'{o}'s planar choreographies, they can be computed to high accuracy using stereographic projection, trigonometric interpolation, and minimization of the action.   

Stability properties of the spherical choreographies have not been discussed. 
In the plane, the only non-circular stable choreography is the figure-eight of Figure~\ref{figure1}. 
We have found numerical evidence that the spherical figure-eight of Figure~\ref{figure5} is stable too.
We have solved the curved $3$-body problem \reff{sphericalnewton}, with initial conditions defined by the reds dots (positions) and the tangents at these dots (velocities) of Figure~\ref{figure5}. 
We ran it for a thousand full orbits, i.e., from $t=0$ to $t=2000\pi$, and the solution did not fall apart.
All the other spherical choreographies presented in this paper fell apart after only a few full orbits.
The systematic approach to study the stability of periodic solutions of dynamical systems is to compute the eigenvalues of the derivatives of the associated Poincar\'{e} maps.
We are currently working on a different algorithm, based on the singular value decomposition of the operator which governs the first variational equation of~\reff{sphericalnewton}, to compute these eigenvalues.
Details will be reported elsewhere.

\section*{Acknowledgements} 

We thank Coralia Cartis and Jared Aurentz for helpful suggestions about quasi-Newton and Newton methods, and Alain Chenciner and Carles Sim\'{o} for giving us details about the computation of planar choreographies.
We are grateful to the reviewers for their comments.
The first author is much indebted to supervisor Nick Trefethen for his continual support and encouragement.

\section*{Appendix A. Closed-form expressions for the gradient and the Hessian in the plane}

Let $N$ be an odd number, and let $p_N(t)$ be the trigonometric interpolant of $q(t)$ at $N$ equispaced points on $[0,2\pi)$ defined by \reff{trigfun}--\reff{trigcoeffs}.
We can decompose the action \reff{action3} into the sum of two terms $A_K$ and $A_U$ with $q(t)$ and $q'(t)$ approximated by $p_N(t)$ and $p_N'(t)$,
\begin{equation}
A_K = \frac{n}{2}\int_0^{2\pi} \big\vert p_N'(t)+ i\omega p_N(t)\big\vert^2 dt, \quad A_U = \frac{n}{2}\sum_{j=1}^{n-1} \int_0^{2\pi} \Big\vert p_N(t) - p_N\Big(t + \frac{2\pi j}{n}\Big) \Big\vert^{-1}dt.
\end{equation}

\noindent The two terms $A_K$ and $A_U$ depend on the $2N$ variables $\{u_k,v_k\}$, $\vert k\vert\leq(N-1)/2$, where $c_k=u_k + iv_k$ are the Fourier coefficients \reff{trigcoeffs} of $p_N(t)$.
Let $\nabla$ denote the gradient with respect to the $u_k$'s and $v_k$'s, that is $\nabla=(\nabla_{u}, \nabla_{v})^T$, $\nabla_{u} =(\partial /\partial u_k)^T$, $\nabla_{v}=(\partial /\partial v_k)^T$, 
$\vert k\vert\leq(N-1)/2$. We wish to derive closed-form expressions for $\nabla A_K$ and $\nabla A_U$.
Consider first $A_K$, with
\begin{equation}
A_K(u_k,v_k) = \pi n \sumodd \vert (k+\omega)c_k\vert^2 = \pi n \sumodd (k+\omega)^2 (u_k^2 + v_k^2),
\end{equation}

\noindent since $p_N'(t)+i\omega p_N(t)$ has Fourier coefficients $\{ikc_k + i\omega c_k\}$, and using Parseval's identity. It leads to
\begin{equation}
\dsp\quad \frac{\partial A_K}{\partial u_k} = 2\pi n(k+\omega)^2u_k, \quad \frac{\partial A_K}{\partial v_k} = 2\pi n(k+\omega)^2v_k, \quad\vert k\vert\leq\frac{N-1}{2}.
\label{dAK}
\end{equation}

Consider now $A_U$, with
\begin{equation}
A_U(u_k,v_k) = \frac{n}{2}\sum_{j=1}^{n-1} \int_0^{2\pi} \frac{dt}{\sqrt{f_j(u_k,v_k,t)}}, \quad \sqrt{f_j(u_k,v_k,t)} = \Big\vert p_N(t) - p_N\Big(t + \frac{2\pi j}{n}\Big) \Big\vert.
\end{equation}

\noindent Expanding $p_N(t)$ and $p_N(t+2\pi j/n)$ and regrouping real and imaginary parts lead to
\begin{equation}
\dsp f_j(u_k,v_k,t) = \dsp\Bigg(\sumodd a_{k,j}(t)u_k + b_{k,j}(t)v_k\Bigg)^2 + \Bigg(\sumodd  a_{k,j}(t)v_k - b_{k,j}(t)u_k\Bigg)^2
\label{fj}
\end{equation}

\noindent with
\begin{equation}
\begin{array}{l}
\dsp a_{k,j}(t) = \big[1-\cos(2\pi jk/n)\big]\cos(kt) + \sin(2\pi jk/n)\sin(kt), \\\\
\dsp b_{k,j}(t) = \big[-1+\cos(2\pi jk/n)\big]\sin(kt) + \sin(2\pi jk/n)\cos(kt),
\label{fjcoeffs}
\end{array}
\end{equation}

\noindent for $\vert k\vert\leq(N-1)/2$ and $1\leq j\leq n-1$. The partial derivatives of $A_U$ with respect to the $u_k$'s and $v_k$'s can then be computed with the chain rule,
\begin{equation}
\dsp \nabla A_U = \frac{n}{2}\sum_{j=1}^{n-1} \int_0^{2\pi}\nabla\bigg(\frac{1}{\sqrt{f_j}}\bigg)dt = -\frac{n}{4}\sum_{j=1}^{n-1} \int_0^{2\pi}\frac{\nabla f_j}{f_j^{3/2}}dt, 
\label{dAU}
\end{equation}

\noindent with
\begin{equation}
\dsp\frac{\partial f_j}{\partial u_k} = 2\sumoddl \Bigg(\Big[a_{l,j}(t)a_{k,j}(t) + b_{l,j}(t)b_{k,j}(t)\Big]u_l + 
\Big[b_{l,j}(t)a_{k,j}(t) - a_{l,j}(t)b_{k,j}(t)\Big]v_l\Bigg),
\label{dfjduk}
\end{equation}

\noindent and
\begin{equation}
\dsp\frac{\partial f_j}{\partial v_k} =  2\sumoddl \Bigg(\Big[a_{l,j}(t)b_{k,j}(t) - b_{l,j}(t)a_{k,j}(t)\Big]u_l + 
\Big[b_{l,j}(t)b_{k,j}(t) + a_{l,j}(t)a_{k,j}(t)\Big]v_l\Bigg).
\label{dfjdvk}
\end{equation}

Let us now derive the formula for the exact Hessian matrix $H$,
\begin{equation}
H = \begin{bmatrix}
\dsp \frac{\partial^2 A}{\partial u_l \partial u_k} & \dsp \frac{\partial^2 A}{\partial u_l \partial v_k}\\\\
\dsp \frac{\partial^2 A}{\partial v_l \partial u_k} & \dsp \frac{\partial^2 A}{\partial v_l \partial v_k}
\end{bmatrix}.
\label{H}
\end{equation}

\noindent $H$ is a $2N\times 2N$ matrix, and each block is $N\times N$. Note that the $\partial^2 A/\partial v_l \partial u_k$ block is the transpose of the $\partial^2 A/\partial u_l \partial v_k$ block, so we are going to derive
formulas for the  $\partial^2 A/\partial u_l \partial u_k$,  $\partial^2 A/\partial u_l \partial v_k$, and $\partial^2 A/\partial v_l \partial v_k$ derivatives only.
It is clear from \reff{dAK} that 
\begin{equation}
\dsp\frac{\partial^2 A_K}{\partial u_l \partial u_k} = \frac{\partial^2 A_K}{\partial v_l \partial v_k} = 2\pi nk^2\delta_{kl}, \quad
\dsp\frac{\partial^2 A_K}{\partial u_l \partial v_k} = 0, \quad
\vert l\vert,\vert k\vert\leq\frac{N-1}{2},
\label{d2AKdvlduk}
\end{equation}

\noindent where $\delta_{kl}$ is the Kronecker delta. 

The second derivatives of $A_U$ with respect to the $u_k$'s and $v_k$'s can be obtained by differentiating \reff{dAU} one more time, e.g.,
\begin{equation}
\dsp \frac{\partial^2 A_U}{\partial u_l\partial u_k} = -\frac{n}{4}\sum_{j=1}^{n-1} \int_0^{2\pi}\frac{\frac{\partial^2 f_j}{\partial u_l \partial u_k} f_j - \frac{3}{2}\frac{\partial f_j}{\partial u_l}\frac{\partial f_j}{\partial u_k}}{f_j^{5/2}}dt, 
\end{equation}

\noindent with 
\begin{equation}
\dsp \frac{\partial^2 f_j}{\partial u_l \partial u_k} = 2(a_{l,j} a_{k,j} + b_{l,j} b_{k,j}).
\end{equation}

\noindent There are similar formulas for the other derivatives with
\begin{equation}
\dsp \frac{\partial^2 f_j}{\partial v_l \partial v_k} = \frac{\partial^2 f_j}{\partial u_l \partial u_k}, \quad \frac{\partial^2 f_j}{\partial u_l \partial v_k} = 2(a_{l,j} b_{k,j} - b_{l,j} a_{k,j}).
\end{equation}

Note that all the second derivatives involving the real and imaginary parts $u_0$ and $v_0$ of the constant terms $c_0$ are zeros, i.e., 
\begin{equation}
\dsp\frac{\partial^2 A_K}{\partial u_l \partial u_0} = \frac{\partial^2 A_K}{\partial v_l \partial v_0} 
= \frac{\partial^2 A_K}{\partial u_l \partial v_0} = 0,
\label{d2AKd0}
\end{equation}

\noindent and 
\begin{equation}
\dsp\frac{\partial^2 A_U}{\partial u_l \partial u_0} = \frac{\partial^2 A_U}{\partial v_l \partial v_0}
= \frac{\partial^2 A_U}{\partial u_l \partial v_0} = 0.
\label{d2AUd0}
\end{equation}

\noindent To prove \reff{d2AKd0}, take $k=0$ in \reff{d2AKdvlduk}, and to prove \reff{d2AUd0}, note that $ a_{0,j} = b_{0,j} = 0$.
As a consequence, when using Newton's method with the exact Hessian \reff{H}, one needs to get rid of these derivatives to make the matrix nonsingular, that is, 
eliminate the lines and columns that correspond to the constant term. Similarly, one needs to get rid of  $u_0$ and $v_0$ in the vector of optimization variables.

\section*{Appendix B. Closed-form expressions for the gradient and the Hessian on the sphere}

Again, let $N$ be an odd number, and let $p_N(t)$ be the trigonometric interpolant of $q(t)$ at $N$ equispaced points on $[0,2\pi)$ defined by \reff{trigfun}--\reff{trigcoeffs}.
We can decompose the action \reff{sphericalaction2bis} into the sum of two terms $A_{K}$ and $A_U$.
The first term comes from the kinetic energy,
\begin{equation}
\dsp A_{K} = \frac{n}{2}\int_0^{2\pi}\bigg(\frac{2R^2\vert p_N'(t) + i\omega p_N(t)\vert}{R^2+\vert p_N(t)\vert^2}\bigg)^2dt,
\end{equation}

\noindent while the second term comes from the potential energy, 
\begin{equation}
\dsp A_U=\frac{n}{2R}\sum_{j=1}^{n-1}\int_0^{2\pi}A_U^{j}(t)dt
\end{equation}

\noindent with
\begin{equation}
\dsp A_U^{j}(t) = \frac{2R^2 - D_j(t)^2}{D_j(t)\sqrt{4R^2-D_j(t)^2}}, \quad 1\leq j\leq n-1,
\end{equation}

\noindent with $D_j(t)=d\big(p_N(t),p_N\big(t+\frac{2\pi j}{n}\big)\big)$.

Again, let $\nabla$ denote the gradient with respect to the $u_k$'s and $v_k$'s, that is $\nabla=(\nabla_{u}, \nabla_{v})^T$, $\nabla_{u} =(\partial /\partial u_k)^T$, $\nabla_{v}=(\partial /\partial v_k)^T$, 
$\vert k\vert\leq(N-1)/2$. Let us first derive the closed-form expression for $\nabla A_{K}$. A straightforward calculation leads to
\begin{equation}
\dsp \nabla A_{K} = \frac{n}{2}\int_0^{2\pi}\nabla\bigg[\bigg(\frac{2R^2\vert p_N'(t) + i\omega p_N(t)\vert}{R^2+\vert p_N(t)\vert^2}\bigg)^2\bigg]dt
= 2nR^4\int_0^{2\pi}\frac{h\nabla g - g\nabla h}{h^2}dt,
\label{dAKsphere}
\end{equation}

\noindent with 
\begin{equation}
\dsp g(u_k,v_k,t) = \vert p_N'(t) + i\omega p_N(t)\vert^2, \quad \dsp h(u_k,v_k,t) = \big(R^2 + \vert p_N(t)\vert^2\big)^2.
\end{equation}

\noindent The functions $g$ and $h$ are given by
\begin{equation}
\dsp g = \Bigg[\sumodd (k+\omega)\Big(u_k\sin(kt) + v_k\cos(kt)\Big)\Bigg]^2 + \Bigg[\sumodd (k+\omega)\Big(u_k\cos(kt) - v_k\sin(kt)\Big)\Bigg]^2,
\end{equation}

\noindent and
\begin{equation}
\dsp h = \Bigg[R^2 + \Bigg(\sumodd u_k\sin(kt) + v_k\cos(kt)\Bigg)^2 + \Bigg(\sumodd u_k\cos(kt) - v_k\sin(kt)\Bigg)^2\Bigg]^2.
\end{equation}

\noindent Their partial derivatives are given by the formulas
\begin{equation}
\begin{array}{l}
\dsp \frac{\partial g}{\partial u_k}  = 2(k+\omega)\dsp\sumoddl 
(l+\omega)\Big(
u_l\cos\big((k-l)t\big)
+ v_l\sin\big((k-l)t\big)
\Big), \\
\dsp \frac{\partial g}{\partial v_k}  = 2(k+\omega)\dsp\sumoddl 
(l+\omega)\Big(
u_l\sin\big((l-k)t\big)
+ v_l\cos\big((l-k)t\big)
\Big),
\end{array}
\end{equation}

\noindent and
\begin{equation}
\begin{array}{l}
\dsp \frac{\partial h}{\partial u_k}  = 4\sqrt{h}\dsp\sumoddl 
u_l\cos\big((k-l)t\big)
+ v_l\sin\big((k-l)t\big), \\
\dsp \frac{\partial h}{\partial v_k}  = 4\sqrt{h}\dsp\sumoddl 
u_l\sin\big((l-k)t\big)
+ v_l\cos\big((l-k)t\big).
\end{array}
\end{equation}

Let us now derive the closed-form expression for $\nabla A_U$, 
\begin{equation}
\dsp \nabla A_U = \frac{n}{2R}\sum_{j=1}^{n-1}\int_0^{2\pi}\nabla A_U^{j}(t)dt, 
\end{equation}

\noindent with 
\begin{equation}
\dsp \nabla A_U^{j}(t) = \frac{-8R^4\nabla D_j(t)}{D_j(t)^2\big[4R^2-D_j(t)^2\big]^{3/2}}, \quad 1\leq j\leq n-1.
\label{dAUjsphere}
\end{equation}

\noindent Let us write
\begin{equation}
D_j(t) = \frac{2R^2\big\vert p_N(t) - p_N(t + \frac{2\pi j}{n}) \big\vert}{\sqrt{\Big(R^2+\big\vert p_N(t)\big\vert^2\Big)\Big(R^2+\big\vert p_N(t + \frac{2\pi j}{n})\big\vert^2\Big)}} = 2R^2\frac{\sqrt{f_j}}{r_0 r_j},
\end{equation}

\noindent with $\sqrt{f_{j}(u_k,v_k,t)} = \vert p_N(t) - p_N(t+2\pi j/n) \vert$ as defined in \reff{fj}--\reff{fjcoeffs}, and $r_j$ defined by
\begin{equation}
r_j(u_k,v_k,t) = \sqrt{R^2 + \Big\vert p_N\Big(t+\frac{2\pi j}{n}\Big)\Big\vert^2},
\end{equation}

\noindent that is,
\begin{equation}
\dsp r_j = \dsp\sqrt{R^2 + \Bigg( \sumodd c_{k,j}(t)u_k + d_{k,j}(t)v_k\Bigg)^2
+\Bigg( \sumodd  c_{k,j}(t)v_k - d_{k,j}(t)u_k \Bigg)^2}, 
\end{equation}

\noindent with
\begin{equation}
\begin{array}{l}
c_{k,j}(t) = \cos(2\pi k j/n)\cos(kt) - \sin(2\pi k j/n)\sin(kt),\\\\
d_{k,j}(t) = -\cos(2\pi k j/n)\sin(kt) - \sin(2\pi k j/n)\cos(kt),
\end{array}
\end{equation}

\noindent for $\vert k\vert\leq(N-1)/2$ and $0\leq j\leq n-1$. It leads to
\begin{equation}
\dsp \nabla D_j(t) = 2R^2\frac{r_0 r_j \nabla f_{j}/2 - f_j\nabla(r_0r_j)}{r_0^2r_j^2\sqrt{f_j}}, \quad 1\leq j\leq n-1.
\end{equation}

\noindent The derivatives of $f_j$ with respect to the $u_k$'s and $v_k$'s are given by \reff{dfjduk}--\reff{dfjdvk}, while the derivatives of $r_j$ are given by
\begin{equation}
\begin{array}{l}
\hspace{-.2cm}\dsp \frac{\partial r_j}{\partial u_k} = \frac{1}{r_j}\Bigg[c_{k,j}(t)\Bigg(\sumoddl c_{l,j}(t)u_l + d_{l,j}(t)v_l\Bigg) 
-d_{k,j}(t)\Bigg(\sumoddl c_{l,j}(t)v_l - d_{l,j}(t)u_l \Bigg)\Bigg],\\
\hspace{-.2cm}\dsp \frac{\partial r_j}{\partial v_k}  = \frac{1}{r_j}\Bigg[d_{k,j}(t)\Bigg(\sumoddl c_{l,j}(t)u_l + d_{l,j}(t)v_l\Bigg)
+ c_{k,j}(t)\Bigg(\sumoddl c_{l,j}(t)v_l - d_{l,j}(t)u_l \Bigg)\Bigg],
\end{array}
\end{equation}

\noindent for $\vert k\vert\leq(N-1)/2$ and $0\leq j\leq n-1$.

Let us now derive the formula for the exact Hessian $H$ on the sphere. Differentiating \reff{dAKsphere} gives the second derivatives of $A_K$ with respect to the $u_k$'s and $v_k$'s, e.g.,
\begin{equation}
\dsp\frac{\partial^2 A_K}{\partial u_l \partial u_k} = 2nR^4\int_0^{2\pi}\frac{\frac{\partial p_k}{\partial u_l}h - 2\frac{\partial h}{\partial u_l}p_k}{h^3}dt,
\quad p_k = h\frac{\partial g}{\partial u_k} - g\frac{\partial h}{\partial u_k},
\end{equation}

\noindent with
\begin{equation}
\frac{\partial p_k}{\partial u_l} = \frac{\partial h}{\partial u_l}\frac{\partial g}{\partial u_k} + h\frac{\partial^2 g}{\partial u_l \partial u_k} - \frac{\partial g}{\partial u_l}\frac{\partial h}{\partial u_k}
-  g\frac{\partial^2 h}{\partial u_l \partial u_k},
\end{equation}

\noindent and
\begin{equation}
\frac{\partial^2 g}{\partial u_l \partial u_k} = 2(k+\omega)(l+\omega)\cos\big((k-l)t\big),
\quad \frac{\partial^2 h}{\partial u_ l \partial u_k} = \frac{1}{2h}\frac{\partial h}{\partial u_l} + 4\sqrt{h}\cos\big((k-l)t\big).
\end{equation}

\noindent There are similar formulas for the other derivatives, with 
\begin{equation}
\frac{\partial^2 g}{\partial v_l \partial v_k} = \frac{\partial^2 g}{\partial u_l \partial u_k},
\quad \frac{\partial^2 g}{\partial u_l \partial v_k} = 2(k+\omega)(l+\omega)\sin\big((l-k)t\big),
\end{equation}

\noindent and 
\begin{equation}
\frac{\partial^2 h}{\partial v_ l \partial v_k} = \frac{1}{2h}\frac{\partial h}{\partial v_l} + 4\sqrt{h}\cos\big((k-l)t\big),
\quad \frac{\partial^2 h}{\partial u_ l \partial v_k} = \frac{1}{2h}\frac{\partial h}{\partial u_l} + 4\sqrt{h}\sin\big((l-k)t\big).
\end{equation}

The second derivatives of $A_U$ can be obtained by differentiating \reff{dAUjsphere}, e.g.,
\begin{equation}
\frac{\partial^2 A_U^{j}}{\partial u_l \partial u_k} = \frac{16R^4\frac{\partial D_j}{\partial u_l}\frac{\partial D_j}{\partial u_k}}{D_j^3\big[4R^2-D_j^2\big]^{3/2}}
- \frac{24R^4\frac{\partial D_j}{\partial u_l}\frac{\partial D_j}{\partial u_k}}{D_j\big[4R^2-D_j^2\big]^{5/2}}
- \frac{8R^4\frac{\partial^2 D_j}{\partial u_l \partial u_k}}{D_j^2\big[4R^2-D_j^2\big]^{3/2}},
\end{equation}

\noindent with 
\begin{equation}
\frac{\partial^2 D_j}{\partial u_l \partial u_k} = 2R^2\frac{r_0r_jf_j\frac{\partial q_k}{\partial u_l} 
- q_k\big(2r_0f_j\frac{\partial r_j}{\partial u_l} + 2r_jf_j\frac{\partial r_0}{\partial u_l} + \frac{1}{2}r_0r_j\frac{\partial f_j}{\partial u_l}\big)}{r_0^3r_j^3f_j^{3/2}},
\label{d2Djdulduk}
\end{equation}

\noindent and
\begin{equation}
q_k = \frac{1}{2}r_0r_j\frac{\partial f_j}{\partial u_k} - f_jr_j\frac{\partial r_0}{\partial u_k} - f_jr_0\frac{\partial r_j}{\partial u_k}.
\label{qk}
\end{equation}

\noindent The second derivatives \reff{d2Djdulduk} involve the derivatives of \reff{qk} with respect of the $u_l$'s---the chain rule leads to nine terms, and only two of them 
are new, $\partial^2 r_0/\partial u_l\partial u_k$ and $\partial^2 r_j/\partial u_l\partial u_k$. These are given by
\begin{equation}
\frac{\partial^2 r_j}{\partial u_l\partial u_k} = \frac{c_{k,j}c_{l,j} + d_{k,j}d_{l,j} - \frac{\partial r_j}{\partial u_l}\frac{\partial r_j}{\partial u_k}}{r_j},
\quad 0\leq j\leq n-1.
\end{equation}

\noindent Similar formulas can be derived for the other derivatives, with 
\begin{equation}
\frac{\partial^2 r_j}{\partial v_l\partial v_k} = \frac{c_{k,j}c_{l,j} + d_{k,j}d_{l,j} - \frac{\partial r_j}{\partial v_l}\frac{\partial r_j}{\partial v_k}}{r_j},
\quad \frac{\partial^2 r_j}{\partial u_l\partial v_k} = \frac{c_{l,j}d_{k,j} - c_{k,j}d_{l,j} - \frac{\partial r_j}{\partial u_l}\frac{\partial r_j}{\partial u_k}}{r_j}.
\end{equation}

Note that, as in the plane, all the second derivatives involving the real and imaginary parts $u_0$ and $v_0$ of the constant terms $c_0$ are zeros, 
so when using Newton's method with the exact Hessian, one needs to get rid of these derivatives to make the matrix nonsingular.

\bibliographystyle{siam}
\bibliography{Montanelli2015}

\begin{thebibliography}{10}

\bibitem{borisov2004}
{\sc A.~V. Borisov, I.~S. Mamaev, and A.~A. Kilin}, {\em Two-body problem on a
  sphere: reduction, stochasticity, periodic orbits}, Regular and Chaotic
  Dynamics, 9 (2004), pp.~265--279.

\bibitem{carinena2005}
{\sc J.~F. Cari\~{n}ena and M.~F. Ra\~{n}ada}, {\em Central potential on spaces
  of constant curvature: the {K}epler problen on the two-dimensional sphere
  $\mathbb{S}^2$ and the hyperbolic plane $\mathbb{H}^2$}, Journal of
  Mathematical Physics, 46 (2005), p.~052702.

\bibitem{chen2008}
{\sc K.-C. Chen}, {\em Existence and minimizing properties of retrograde orbits
  to the three-body problem with various choices of masses}, Annals of
  Mathematics, 167 (2008), pp.~325--348.

\bibitem{chenciner2000a}
{\sc A.~Chenciner and R.~Montgomery}, {\em A remarkable periodic solution of
  the three-body problem in the case of equal masses}, Annals of Mathematics,
  152 (2000), pp.~881--901.

\bibitem{diacu2012c}
{\sc F.~Diacu}, {\em Relative equilibria of the curved {N}-body problem},
  Springer, 2012.

\bibitem{diacu2013c}
{\sc F.~Diacu and S.~Kordlou}, {\em Rotopulsators of the curved {N}-body
  problem}, Journal of Differential Equations, 255 (2013), pp.~2709--2750.

\bibitem{diacu2013b}
{\sc F.~Diacu, R.~Mart\'{\i}nez, E.~P\'{e}rez-Chavala, and C.~Sim\'{o}}, {\em
  On the stability of tetrahedral relative equilibria in the positively curved
  4-body problem}, Physica D, 256-257 (2013), pp.~21--35.

\bibitem{diacu2011}
{\sc F.~Diacu and E.~P\'{e}rez-Chavala}, {\em Homographic solutions of the
  curved 3-body problem}, Journal of Differential Equations, 250 (2011),
  pp.~340--366.

\bibitem{diacu2012a}
{\sc F.~Diacu, E.~P\'{e}rez-Chavala, and M.~Santoprete}, {\em The n-body
  problem in spaces of constant curvature. {P}art {I}: relative equilibria},
  Journal of Nonlinear Science, 22 (2012), pp.~247--266.

\bibitem{diacu2012b}
\leavevmode\vrule height 2pt depth -1.6pt width 23pt, {\em The n-body problem
  in spaces of constant curvature. {P}art {II}: singularities}, Journal of
  Nonlinear Science, 22 (2012), pp.~267--275.

\bibitem{driscoll2008}
{\sc T.~A. Driscoll, F.~Bornemann, and L.~N. Trefethen}, {\em {The chebop
  system for automatic solution of differential equations}}, BIT Numerical
  Mathematics, 48 (2008), pp.~701--723.

\bibitem{chebfun}
{\sc T.~A. Driscoll, N.~Hale, and L.~N. Trefethen}, eds., {\em {C}hebfun
  {G}uide}, Pafnuty Publications, 2014.

\bibitem{lagrange1772}
{\sc J.-L. Lagrange}, {\em Essai sur le probl\`{e}me des trois corps}, in Prix
  de l'Acad\'{e}mie Royale des Sciences, vol.~IX, 1772, pp.~229--332.

\bibitem{maupertuis1744}
{\sc P.~L. Maupertuis}, {\em Accord de diff\'{e}rentes loix de la nature qui
  avoient jusqu'ici paru incompatibles}, M\'{e}moires de l'Acad\'{e}mie Royale
  des Sciences,  (1744), pp.~417--426.

\bibitem{moore1993}
{\sc C.~Moore}, {\em Braids in classical dynamics}, Physical Review Letters, 70
  (1993), pp.~3675--3679.

\bibitem{nocedal2006}
{\sc J.~Nocedal and S.~J. Wright}, {\em Numerical Optimization}, Springer,
  second~ed., 2006.

\bibitem{perez2012}
{\sc E.~P\'{e}rez-Chavala and J.~G. Reyes-Victoria}, {\em An intrinsic approach
  in the curved n-body problem. {T}he positive curvature case}, Transactions of
  the American Mathematical Society, 364 (2012), pp.~3805--3827.

\bibitem{poincare1892}
{\sc H.~Poincar\'{e}}, {\em Les M\'{e}thodes Nouvelles de la M\'{e}canique
  C\'{e}leste}, vol.~I: Solutions p\'{e}riodiques, Non-existence des
  int\'{e}grales uniformes, Solutions asymptotiques, Gauthier-Villars et Fils,
  1892.

\bibitem{poincare1896}
\leavevmode\vrule height 2pt depth -1.6pt width 23pt, {\em Sur les solutions
  p\'{e}riodiques et le principe de moindre action}, Comptes Rendus de
  l'Acad\'{e}mie des Sciences, 123 (1896).

\bibitem{shanno1970}
{\sc D.~Shanno}, {\em Conditioning of quasi-{N}ewton methods for function
  minimization}, Mathematics of Computation, 24 (1970), pp.~647--656.

\bibitem{simo2001}
{\sc C.~Sim\'{o}}, {\em New families of solutions in {N}-body problems}, in
  Proceedings of the Third European Congress of Mathematics, Birkh\"{a}user
  Verlag, 2001.

\bibitem{trefethen2014}
{\sc L.~N. Trefethen and J.~A.~C. Weideman}, {\em The exponentially convergent
  trapezoidal rule}, SIAM Review, 56 (2014), pp.~385--458.

\bibitem{wright2015}
{\sc G.~B. Wright, M.~Javed, H.~Montanelli, and L.~N. Trefethen}, {\em
  Extension of {C}hebfun to periodic functions}, SIAM Journal on Scientific
  Computing, to appear.

\end{thebibliography}

\end{document}